\documentclass[runningheads]{llncs}
\usepackage[T1]{fontenc}
\usepackage{graphicx}
\usepackage[numbers,sort&compress]{natbib}
\usepackage{hyperref}
\usepackage{xcolor}
\usepackage{booktabs}
\usepackage{listings}
\usepackage{amsmath}
\usepackage{algorithm}
\usepackage{algorithmic}
\usepackage{amssymb}
\usepackage{subcaption}
\usepackage{makecell}

\definecolor{mygreen}{rgb}{0,0.6,0}
\definecolor{mygray}{rgb}{0.5,0.5,0.5}
\definecolor{mymauve}{rgb}{0.58,0,0.82}

\lstdefinestyle{Pythonstyle}{
  language=Python,
  basicstyle=\ttfamily\footnotesize,
  keywordstyle=\color{violet}\bfseries,
  commentstyle=\color{mygreen}\textit,
  stringstyle=\color{mymauve},
  showstringspaces=false,
  breaklines=true,
  frame=tb,
  numbers=left,
  numberstyle=\tiny\color{mygray},
}

\definecolor{revtwocolor}{rgb}{0.0, 0.48, 0.65}

\renewcommand\refname{References}

\newcommand{\del}[1]{{\iffalse {#1} \fi}}

\newcommand{\revone}[1]{{\color{blue}{#1}}}
\newcommand{\revtwo}[1]{{\color{revtwocolor}{#1}}}

\renewcommand{\revone}[1]{{\color{black}{#1}}}
\renewcommand{\revtwo}[1]{{\color{black}{#1}}}
\newcommand{\rev}[1]{{\color{black}{#1}}}

\begin{document}
\title{Addressing Methodological \rev{Sensitivity} in MCDM with a Systematic Pipeline Approach to Data Transformation Sensitivity Analysis}
%
%
\author{Juan B. Cabral\inst{1,2,3}\orcidID{0000-0002-7351-0680} \and
Alvaro Roy Schachner\inst{3,4}\orcidID{0009-0008-5460-8720}}
%
\authorrunning{Cabral \& Schachner, 2026}
\titlerunning{\rev{Methodological Sensitivity in MCDM}}
%
\institute{Grupo de Innovación y Desarrollo Tecnológico, Comisión Nacional de Actividades Espaciales (GVT-CONAE), Córdoba, Argentina \and
Consejo Nacional de Investigaciones Científicas y Técnicas (CONICET), Córdoba, Argentina \and
Facultad de Matemática Astronomía y Física, Universidad Nacional de Córdoba (FAMAF-UNC), Córdoba, Argentina \and
Instituto de Astronomía Téorica y Experimental (IATE-CONICET), Córdoba, Argentina \\
\email{jbcabral@unc.edu.ar - alvaro.schachner@mi.unc.edu.ar}}
\maketitle              

\begin{abstract}
Multicriteria decision-making methods exhibit critical dependence on the choice of normalization techniques, where different selections can alter 20–40\% of the final rankings. Current practice is characterized by the ad-hoc selection of methods without systematic robustness evaluation. We present a framework that addresses this methodological \rev{sensitivity} through automated exploration of the scaling transformation space. The implementation leverages the existing Scikit-Criteria infrastructure to automatically generate all possible methodological combinations and provide robust comparative analysis.
We apply this approach in an evaluation dataset of cryptocurrencies with 6 methodological scenarios, showing a range of correlation between methods, explicitly quantifying the methodological \rev{sensitivity} limits.

\end{abstract}

\section{Introduction}

Multi Criteria Decision Methods (MCDM) have long been affected by
a fundamental limitation: a {critical dependency of the
final results on methodological decisions}; in particular, the
election of the data scaling methods \citep{chen2019effects}.
Empirical evidence shows that different normalization methods
could alter around {20\% to 40\% of the final rankings}
in real applications~\citep{aytekin2021comparative, krishnan2022past}.
However, current practices are often characterized by the
\textit{ad-hoc} selection of scaling methods, based on
methodological precedents, personal preference,
or software availability, without a systematic evaluation of
their impact on the robustness of the results.
Some limitations of prevailing approaches can be outlined as
follows:
1) The MCDM literature lacks \textit{unified conceptual frameworks} for systematically evaluating sensitivity to scaling transformations. Existing studies are fragmented, employ heterogeneous methodologies, and are often limited to pairwise comparisons between specific methods without considering the broader methodological landscape \citep{mukhametzyanov2021specific}.
2) Traditional approaches generally provide point estimates
without quantifying \rev{sensitivity}, which may obscure the inherent
variability introduced by methodological choices; particularly
in contexts where methodological transparency is essential.
3) Available software implementations of MCDM methods typically
require users to specify a single normalization procedure prior to the analysis, which restricts the possibility of
systematically exploring the methodological space \citep{kizielewicz2023pymcdm},
among other limitations.
This work addresses the fundamental research question ``How can the variability in rankings be \revone{automatized}, quantified and communicated through combinations of normalization and aggregation \revone{in decision making}?''
\rev{To overcome these limitations, we propose a combinatorial approach}, in the mathematical sense of systematically exploring all possible combinations from finite sets of methodological choices, that enables exhaustive evaluation of scaling and aggregation method combinations, providing explicit quantification of methodological \rev{sensitivity} bounds.

Tools derived from techniques such as Grid Search and their
implementation in Scikit-Learn \citep{pedregosa2011scikit} illustrate how the field of machine learning has successfully adopted combinatorial pipeline architectures. This design pattern serves as a reference point for our proposed solution and conceptual validation in the MCDM domain, highlighting the feasibility of systematically exploring complex methodological spaces through the automated combination of algorithmic components. Moreover, by leveraging modern parallel computing frameworks, it becomes possible to overcome the computational barriers that have historically limited the practicality of exhaustive combinatorial approaches.

In this work, we introduce a systematic framework for the analysis of MCDM data transformations, which enables the automated exploration of the methodological space and the explicit quantification of its \rev{sensitivity}, built upon the existing infrastructure of \textit{Scikit-Criteria}. In particular, the \textit{RanksComparator} proved especially valuable, as it already facilitates the comparison and  evaluation of different rankings \citep{cabral2025beyond}.

Our contributions are threefold.
1) An algorithmic framework for the automated exploration of the MCDM transformation space through the Cartesian product of modular components (\texttt{SKCCombinatorialPipeline}).
2) A multi-metric \revtwo{sensitivity analysis} approach that provides correlations, covariances, coefficients of determination ($R^2$), and distance-based measures between rankings generated by different methodological configurations, leveraging the existing ranking comparison infrastructure.
3) An empirical validation on a real-world cryptocurrency dataset, yielding interpretable \revtwo{sensitivity} bounds that enable the assessment of the robustness of the resulting conclusions.

\color{black}

\section{Related Work: Normalization Methods and Their Impact on MCDM Rankings}

Although it is intuitive that different aggregation functions may lead to different rankings, empirical evidence has shown that alternative normalization methods can also produce inconsistent results across multiple aggregation functions~\citep{opricovic2004compromise}.
Subsequent studies have quantified this variability, reporting that the choice of normalization method can affect between 20\% and 40\% of final rankings in real-world applications \citep{aytekin2021comparative, krishnan2022past}.

Moreover, a comprehensive survey identified 31 distinct normalization methods used in the MCDM literature, each exhibiting different mathematical properties and behavioral characteristics \citep{jahan2015state}.
In this context, method-specific recommendations have been developed regarding the choice of normalization techniques for particular MCDM methods, based on extensive empirical validations \citep{chen2019effects, krishnan2022past}.

Nevertheless, a critical gap has also been identified in the synthesis of 19 comparative studies: \textit{no study has simultaneously evaluated all normalization methods in combination with all commonly used MCDM techniques}.
This gap directly motivates the contribution of our exhaustive combinatorial approach, which enables such analyses to be conducted in a straightforward and reproducible manner, even in practical settings where methodological evaluations must remain simple and accessible \citep{krishnan2022past}.

\subsection{Software Tools for MCDM}

\autoref{tab:software_comparison} presents a summary on the capabilities of existing MCDM tools.
PyMCDM \citep{kizielewicz2023pymcdm} provides more than 15 MCDM methods but requires manual comparison across methodological configurations.
D-Sight and Visual PROMETHEE are commercial tools specialized in the PROMETHEE family of methods.
DIVIZ/Decision Deck \citep{cinelli2022proper} supports complex workflows through XML-based specifications; however, its development has remained largely inactive since approximately 2019.
The MCDM package for R implements a MetaRanking approach using a fixed set of predefined methods.

\begin{table}
\centering
\caption{Comparison of MCDM software tools.}
\label{tab:software_comparison}
\begin{tabular}{@{}lccc@{}}  
\toprule
\textbf{Tool} & \textbf{Methods} & \textbf{Combinatorial} & \textbf{Status} \\
\midrule
PyMCDM & 15+ & Manual only & Active \\
D-Sight/V-PROMETHEE & PROMETHEE & No & Commercial \\
DIVIZ/Decision Deck & Multiple & XML workflows & Frozen \\
R MCDM package & MetaRanking & Fixed set & Active \\
\midrule  
\textbf{Scikit-Criteria} & \textbf{20+} & \textbf{Automated} & \textbf{Active} \\
\bottomrule
\end{tabular}
\vspace{11pt}
\end{table}

To the best of our knowledge, no existing tool applies a GridSearchCV-like exploration strategy (inspired by Scikit-Learn \citep{pedregosa2011scikit}) to the MCDM methodological space.
This constitutes the distinctive contribution of \texttt{SKCCombinatorialPipeline}, which enables an exhaustive and automated exploration of the Cartesian product of methodological configurations.

\section{Scikit-Criteria Rank Comparison Infrastructure}
\label{section:rank_cmp}

\textit{Scikit-Criteria} has evolved into a comprehensive library providing a powerful set of tools for discrete multi-criteria decision-making methods \citep{cabral2016scikit}.
In this work, we specifically built upon the existing infrastructure: the \texttt{RanksComparator}.

\subsection{Formal Definitions}

In Scikit-Criteria, a ranking result is defined as follows:

\begin{definition}[RankResult]
A \texttt{RankResult} is a tuple
\begin{equation}
RankResult = (M, A, r, E)
\end{equation}
where:
\begin{itemize}
\item $M$: method identifier (string)
\item $A = \{a_1, \ldots, a_m\}$: set of alternatives
\item $r: A \to \mathbb{N}$: ranking function, where $r(a_i)$ represents the position of $a_i$ (lower values = better ranking)
\item $E$: additional metadata dictionary (intermediate calculations, parameters, and other non-structured data for interpreting the ranking)
\end{itemize}
\end{definition}

The \texttt{RankResult} structure represents the output of MCDM methods that produce ordered rankings of alternatives.
It incorporates method metadata, intermediate calculations, and also provides the ability to generate untied rankings through the \texttt{untied\_rank\_} property, which resolves ties by following the original order of alternatives while preserving their relative positions.

A ranks comparator is defined as:

\begin{definition}[RanksComparator]
A \texttt{RanksComparator} is a pair
\begin{equation}
RanksComparator = (\mathcal{R}, E)
\end{equation}
where:
\begin{itemize}
\item $\mathcal{R} = \{(N_1, R_1), \ldots, (N_k, R_k)\}$: set of $k \geq 2$ named rankings
\item $N_i$: unique name of the $i$-th ranking (string)
\item $R_i = (M_i, A, r_i, E_i)$: $i$-th \texttt{RankResult}
\item $E$: additional metadata dictionary for the comparator (extra)
\end{itemize}
such that:
\begin{itemize}
\item $\forall i,j \in \{1,\ldots,k\}: R_i.A = R_j.A$ (all rankings have the same alternatives)
\item $\forall i \neq j: N_i \neq N_j$ (unique ranking names)
\item $k \geq 2$ (at least two ranks to compare)
\end{itemize}
\end{definition}

This construct \texttt{RanksComparator} enables the comparison of multiple rankings over the same set of alternatives \citep{cabral2025beyond}, providing comprehensive \revtwo{sensitivity analysis} through multiple complementary metrics: \texttt{corr()} which allows computing \revone{Spearman, Kendall, Pearson correlation or defining custom metrics} for pairwise rankings, measuring \revone{monotonic} agreement between methodological combinations; \texttt{cov()} calculates pairwise covariance of rankings, quantifying joint variability between methods; \texttt{r2\_score()} computes pairwise coefficient of determination regression score function, measuring how well one ranking can predict another;
\revone{The \texttt{distance()} function supports 20 pairwise distance metrics between rankings, including Hamming (normalized), Euclidean, and Manhattan distances, while also allowing the definition of custom distance functions that provide a geometric interpretation of methodological differences.}
%
Additionally, it offers visualizations via flow diagrams, heatmaps, and boxplots, as well as data conversion to widely adopted Python types such as \texttt{pandas.Data-\allowbreak Frame} \citep{mckinney2011pandas}. This multidimensional characterization enables decision-makers to establish bounds on methodological \revtwo{sensitivity} and assess the robustness of their conclusions.
This existing infrastructure provides the base architecture for our combinatorial framework, where we extend our comparison capabilities toward a systematic exploration over some complete methodological space over \texttt{SKCCombinatorial-\allowbreak Pipeline}.

\section{Combinatorial Pipelines}
\label{section:comb_pipelines}

For several versions of Scikit-Criteria, \texttt{SKCPipeline} allowed composing processing pipelines where each step transforms the decision matrix until reaching a final aggregator that produces a \texttt{RankResult} (see \autoref{code:pipe}). This sequential architecture works with predefined individual steps.

\revone{The pipeline supports various preprocessing transformers and aggregation methods. Among the aggregators, the \textbf{Weighted Sum Model (WSM)} computes a weighted linear combination of normalized criteria values, while \textbf{TOPSIS} (Technique for Order of Preference by Similarity to Ideal Solution) ranks alternatives by their distance to an ideal solution. Scalers include \textbf{SumScaler} \rev{(Sum)}, which normalizes each criterion by the sum of its values; \textbf{VectorScaler} \rev{(Vec)}, which uses the Euclidean norm; and \textbf{MinMaxScaler} \rev{(MM)}, which maps values to the $[0,1]$ interval using minimum and maximum values.}

\begin{figure}[h]
\vspace{11pt}
\begin{lstlisting}[language=Python]
from skcriteria.pipeline import mkpipe  # import pipeline
pipeline = mkpipe(  # create pipeline
   NegateMinimize(),  # negate to maximize all criteria
   FilterGT({'criteria': 300}),  # apply satisficing filter
   FilterNonDominated(),  # remove dominated alternatives
   SumScaler(target="weights"),  # normalize weights proportionally
   MinMaxScaler(target="matrix"),  # normalize matrix by min-max
   TOPSIS()  # apply TOPSIS method
)
pipeline.evaluate(dm)  # evaluate some decision-matrix
\end{lstlisting}
\caption{Modular MCDM pipeline in Scikit-Criteria that shows the systematic composition of data transformations, alternative filters, and normalization methods, prior to the final evaluation.}
\label{code:pipe}
\end{figure}

The new \texttt{SKCCombinatorialPipeline} implements true combinatorial
exploration by systematically generating the Cartesian product of
methodological alternatives at each pipeline step. When a step contains
$n$ alternative methods and the next step contains $m$ alternatives,
the system automatically constructs all $n \times m$ possible
combinations, ensuring exhaustive coverage of the methodological space.

For example, by defining two steps: \texttt{scaler} with two alternatives \texttt{SumScaler} and \texttt{VectorScaler}; and \texttt{agg} with two options \texttt{WSM} and \texttt{TOPSIS}; the system automatically generates four pipelines:

\begin{enumerate}
\item \texttt{SumScaler+WSM}
\item \texttt{SumScaler+TOPSIS}
\item \texttt{VectorScaler+WSM}
\item \texttt{VectorScaler+TOPSIS}.
\end{enumerate}

This architecture {fully leverages} the existing \textit{Scikit-Criteria} infrastructure, particularly \texttt{RanksComparator}, to provide automatic methodological sensitivity analysis through systematic generation and comparison of all possible methodological combinations.
%
\autoref{alg:combinatorial} formally describes the combinatorial evaluation process.

\begin{algorithm}
\begin{algorithmic}[1]
\STATE \textbf{Input:} $dm$ (DecisionMatrix), $steps = [S_1, S_2, \ldots, S_k]$
\STATE \textbf{Output:} RanksComparator
\STATE $combinations \gets S_1 \times S_2 \times \cdots \times S_k$ \COMMENT{Cartesian product}
\STATE $results \gets [ ]$ \COMMENT{Empty list}
\FOR{each $config$ in $combinations$}
    \STATE $pipeline \gets$ SKCPipeline($config$) \COMMENT{Instantiate pipeline}
    \STATE $rank \gets pipeline$.evaluate($dm$) \COMMENT{Parallel evaluation}
    \STATE $name \gets$ generateName($config$) \COMMENT{Unique identifier}
    \STATE $results$.append(($name$, $rank$))
\ENDFOR
\STATE \textbf{return} RanksComparator($results$)
\end{algorithmic}
\caption{SKCCombinatorialPipeline Evaluation.}
\label{alg:combinatorial}
\end{algorithm}

\revone{The \texttt{generateName()} function creates unique identifiers by concatenating the class names of the components in each step (e.g., \texttt{SumScaler\_TOPSIS}). A numbered suffix (e.g., \texttt{\_1}, \texttt{\_2}) is appended only when two different pipeline configurations produce identical string representations, a situation that may arise when user-defined components share the same class name or when the same class is instantiated with different parameters that are not reflected in its string representation.}
The algorithm works through a systematic four-stage process: 1) combinatorial generation by applying the Cartesian product over collections of algorithmic components, 2) pipeline instantiation through automatic construction of \texttt{SKCPipeline} objects with unique identifiers for each generated combination, 3) distributed evaluation via parallel invocation of each constructed pipeline, returning a \texttt{RanksComparator} object, and 4) infrastructure integration through transparent use of the existing API for ranking comparison and analysis.
\subsection{Complexity Analysis}

The practical viability of the combinatorial approach greatly depends on its computational cost.
\revtwo{
If we consider a pipeline with $k$ steps $S_1,S_2,...,S_k$, in a bounded parallel execution context, the temporal complexity is:

\begin{equation}
    T_{parallel} = O \left( \frac{S}{P} \cdot \sum_{i=1}^{k} \frac{1}{|S_i|} \sum_{j=1}^{|S_i|} c_{i,j} \right),
\end{equation}
where $P$ is the number of processors, $S=\prod_{i=1}^{k} |S_i|$ is the total number of paths and $c_{i,j}$ is the cost of alternative $j$ of step $i$.
If we consider the special case of a uniform cost $c$ across all stages, the resulting complexity can be reduced to:
\begin{equation}
    T_{parallel} = O \left( \frac{S \, k \, c}{P} \right),
\end{equation}
for instance, we can use $c=m \times n$ (alternatives $\times$ criterion) with aggregators such as WSM and WPM.
}
Since each pipeline is independent, evaluation is \textit{embarrassingly parallel}, approaching an ideal linear speedup.

For example, given $k=2$ steps with $|S_1|=3$ normalizers and $|S_2|=2$ aggregators, $m=10$ alternatives, $n=5$ criterion and $P=4$ cores:

\begin{equation}
T_{parallel} = O\left(\frac{50 \times 2 \times 6}{4}\right) = O(150)
\end{equation}

Complexity grows exponentially with $k$, but stays manageable for typical methodological spaces ($k \leq 4$ steps, $|S_i| \leq 5$ options per-step).
\revtwo{It is also worth noting that the case study presented in Section~\ref{section:case_study} (6 pipelines on a $9 \times 6$ matrix) represents a deliberately simple, illustrative scenario. Its computational cost is intentionally trivial; the framework's scalability to larger combinatorial spaces is supported by the parallel architecture described above.}

\subsection{Extensibility}

Scikit-Criteria follows a uniform interface.
Transformers (e.g., normalizers) are classes that implement a \texttt{transform(dm) $\to$ \revtwo{DecisionMatrix}} method, and can be used in the preprocessing phase of a pipeline.
Similarly, aggregators are classes that implement an \texttt{evaluate(dm) $\to$ RankResult} method, and can be used in the final phase of a pipeline.

All existing infrastructure within Scikit-Criteria is compatible with this architecture, that is, transformers (filters, inverters, scalers) and aggregators.
This modular approach allows transparent integration of new community-made methods without modifying the core combinatorial framework.


\section{Case Study: Cryptocurrency Evaluation}
\label{section:case_study}

To demonstrate the practical application of our framework, we present a scaling and aggregation function sensitivity analysis using the cryptocurrency evaluation dataset from Van Heerden et al.~\citep{van2021evaluation}.
\revone{Table~\ref{tab:crypto_names} lists the nine cryptocurrencies included in the dataset along with their full names.}
The dataset measures these assets across multiple criteria including market capitalization, volatility, and trading volume.

\revone{
\begin{table}
\centering
\caption{Cryptocurrency ticker codes and corresponding full names used in the dataset.}
\label{tab:crypto_names}
\begin{tabular}{ll|ll}
\toprule
\textbf{Ticker} & \textbf{Full Name} & \textbf{Ticker} & \textbf{Full Name} \\
\midrule
ADA  & Cardano & BNB  & Binance Coin \\
BTC  & Bitcoin & DOGE & Dogecoin \\
ETH  & Ethereum & LINK & Chainlink \\
LTC  & Litecoin & XLM  & Stellar \\
XRP  & Ripple & & \\
\bottomrule
\end{tabular}
\end{table}
}

\begin{figure}
    \centering
    \vspace{5pt}
    \includegraphics[width=1\textwidth]{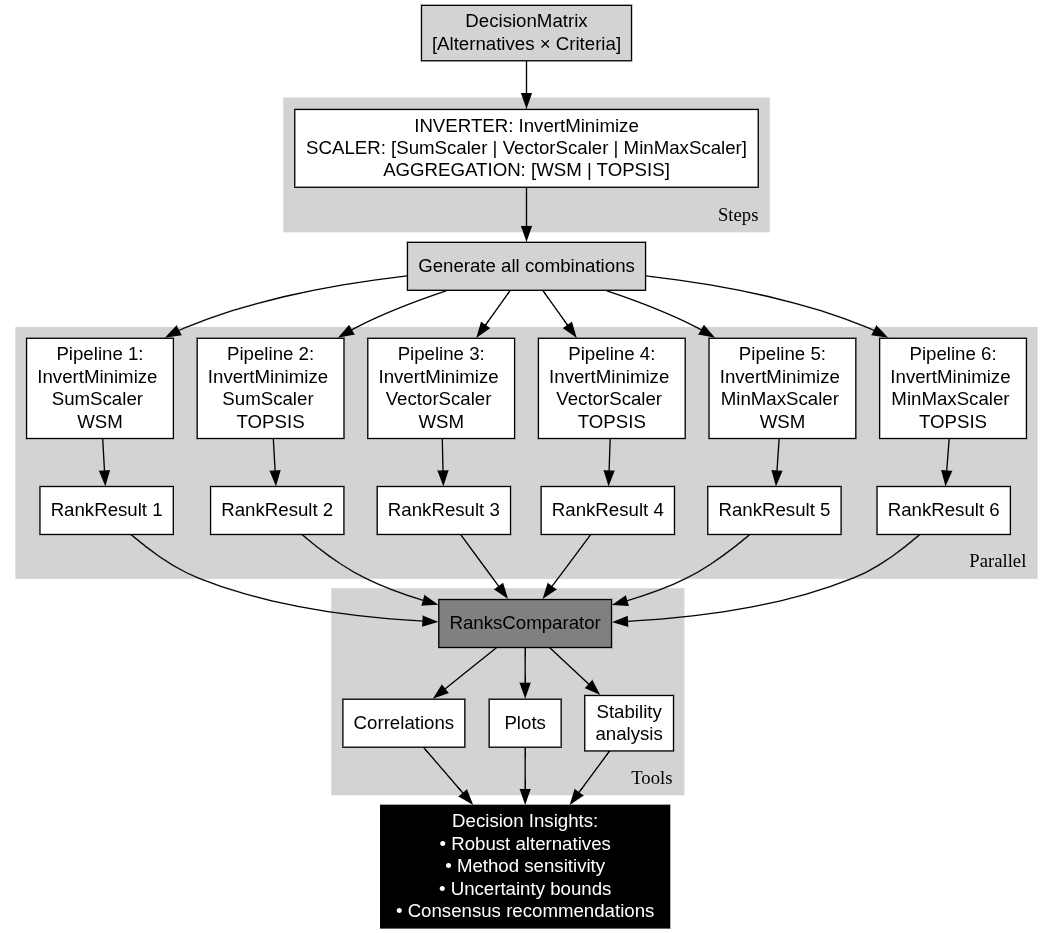}
    \caption{Combinatorial pipeline framework architecture for systematic MCDM sensitivity analysis.}
    \label{fig:pipeline}
\end{figure}

The experiment consists of loading this dataset and creating a combinatorial pipeline that includes \revone{an \textbf{inversion procedure} for minimization criteria (\texttt{InvertMinimize})}, three scaling options (\texttt{SumScaler}, \texttt{VectorScaler}, and \texttt{MinMaxScaler} applied to the data matrix), and two aggregation methods (\texttt{Weig-\allowbreak htedSumModel} and \texttt{TOPSIS}), automatically generating \textit{six} different methodological scenarios (3 scalers $\times$ 2 aggregation methods) that evaluate each combination and visualize the resulting ranking distributions.
\revone{The \texttt{InvertMinimize} transformer converts minimization criteria to maximization by computing the reciprocal of each value ($c' = 1/c$), which \rev{do not preserve proportionality} for strictly positive values. Analysts should be aware that this transformation may not preserve scale characteristics when criteria include zero or negative values (for more details on this transformation, refer to the \texttt{scikit-criteria} documentation \footnote{\url{https://scikit-criteria.quatrope.org/en/latest/api/preprocessing/invert_objectives.html}}).}
This algorithmic idea can be appreciated in more detail in \autoref{fig:pipeline}, its Python code implementation in \autoref{code:comb}.

\begin{figure}[h]
\begin{lstlisting}[language=Python]
dm = skc.datasets.load_van2021evaluation()  # Load dataset
pipeline = mkcombinatorial(  # create combinatorial pipeline
  InvertMinimize(),
  [
    SumScaler(target="matrix"),
    VectorScaler(target="matrix"),
    MinMaxScaler(target="matrix")
  ],
  [WeightedSumModel(), TOPSIS()]
)
rank_comparator = pipeline.evaluate(dm) # Evaluate
rank_comparator.plot()  # plot
rank_comparator.corr("spearman")  # correlation matrix
\end{lstlisting}
\vspace{-2pt}
\caption{Example use of the combinatorial framework.}
\label{code:comb}
\end{figure}

\begin{table}
\centering
\caption{Rankings of cryptocurrencies across six methodological configurations with descriptive statistics.}
\label{tab:ranks}
\begin{tabular}{l|rrrrrr@{\hspace{4pt}}|@{\hspace{4pt}}r@{\hspace{8pt}}rrr}
\toprule
\makecell{\\ \textbf{Alternative}} &
\makecell{Sum \\ WSM} &
\makecell{Sum \\ TOPSIS} &
\makecell{Vec \\ WSM} &
\makecell{Vec \\ TOPSIS} &
\makecell{MM \\ WSM} &
\makecell{MM \\ TOPSIS} &
Median & MAD & min & max \\
\midrule
ADA  & 8 & 8 & 8 & 8 & 8 & 8 & 8.000 & 0.000 & 8 & 8 \\
BNB  & 2 & 2 & 2 & 2 & 1 & 1 & 2.000 & 0.000 & 1 & 2 \\
BTC  & 1 & 1 & 1 & 1 & 2 & 2 & 1.000 & 0.000 & 1 & 2 \\
DOGE & 7 & 5 & 7 & 5 & 6 & 4 & 5.500 & 1.000 & 4 & 7 \\
ETH  & 3 & 3 & 3 & 4 & 3 & 5 & 3.000 & 0.000 & 3 & 5 \\
LINK & 5 & 6 & 5 & 6 & 4 & 6 & 5.500 & 0.500 & 4 & 6 \\
LTC  & 6 & 7 & 6 & 7 & 7 & 7 & 7.000 & 0.000 & 6 & 7 \\
XLM  & 4 & 4 & 4 & 3 & 5 & 3 & 4.000 & 0.500 & 3 & 5 \\
XRP  & 9 & 9 & 9 & 9 & 9 & 9 & 9.000 & 0.000 & 9 & 9 \\
\bottomrule
\end{tabular}
\vspace{11pt}
\end{table}

\color{black}

\begin{figure}
    \centering
    \includegraphics[width=\textwidth]{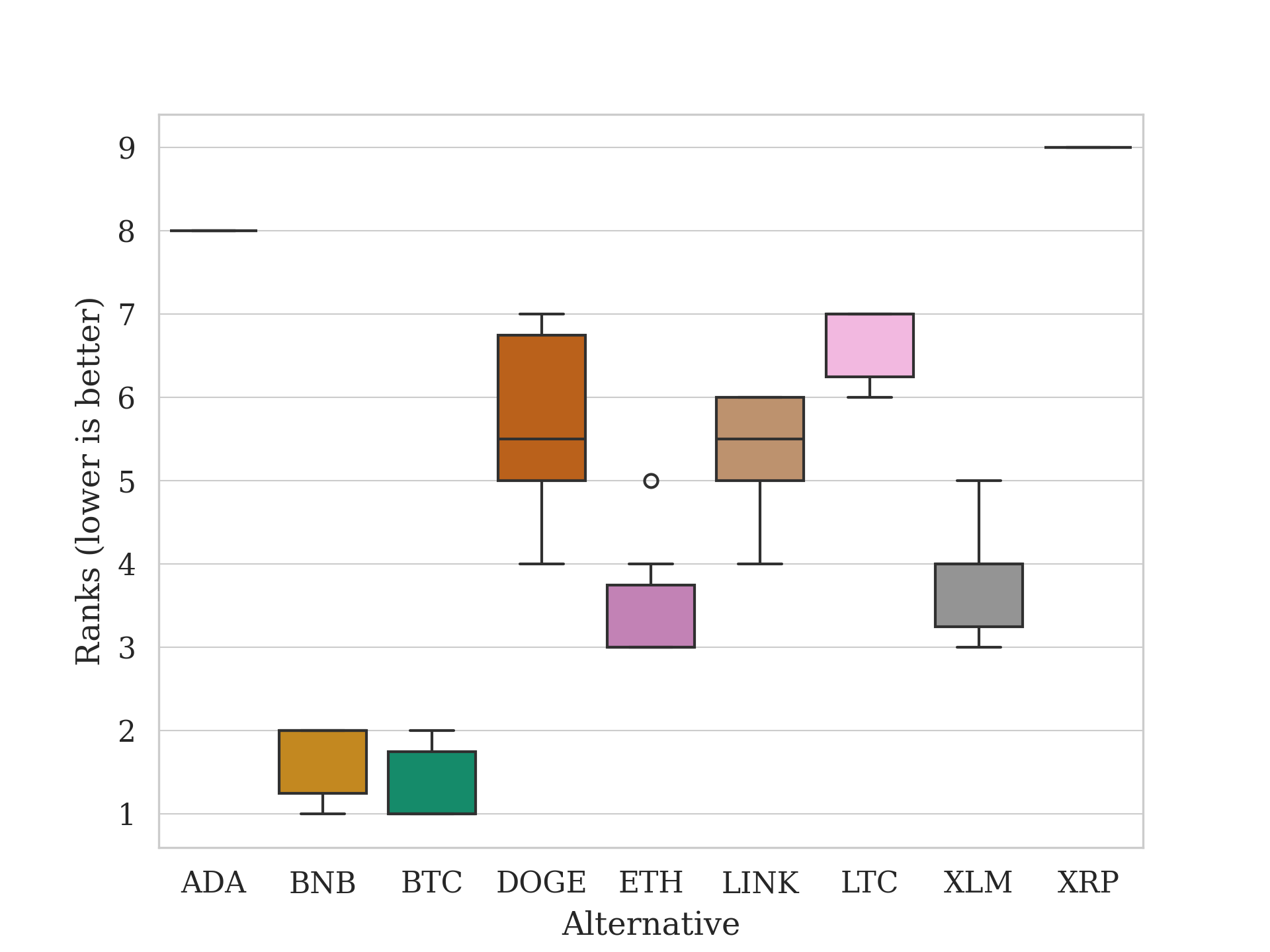}\\
    \vspace{.5em}
    \hrule
    \vspace{.5em}
    \begin{tabular}{l|rrrrrr}
\toprule
\textbf{Method} & \makecell{Sum\\WSM} & \makecell{Sum\\TOPSIS} & \makecell{Vec\\WSM} & \makecell{Vec\\TOPSIS} & \makecell{MM\\WSM} & \makecell{MM\\TOPSIS} \\
\midrule
Sum WSM & 1.000 & 0.950 & 1.000 & 0.933 & 0.950 & 0.850 \\
Sum TOPSIS & 0.950 & 1.000 & 0.950 & 0.983 & 0.933 & 0.933 \\
Vec WSM & 1.000 & 0.950 & 1.000 & 0.933 & 0.950 & 0.850 \\
Vec TOPSIS & 0.933 & 0.983 & 0.933 & 1.000 & 0.900 & 0.967 \\
MM WSM & 0.950 & 0.933 & 0.950 & 0.900 & 1.000 & 0.867 \\
MM TOPSIS & 0.850 & 0.933 & 0.850 & 0.967 & 0.867 & 1.000 \\
\bottomrule
\end{tabular}
    \caption{\revone{Output of the Cryptocurrency Evaluation: (a) Boxplot and (b) Spearman Correlation Matrix.}}
    \label{fig:code_plot}
\end{figure}

From the box plot in \autoref{fig:code_plot}(a)\revone{, several key insights can be observed.}
\revone{BNB and BTC consistently achieve superior positions (ranks 1--2) across all methodological combinations, indicating that these alternatives are robust options regardless of the decision approach used. Conversely, ADA and XRP remain perfectly stable at positions 8 and 9 respectively across all six configurations (MAD $= 0$), indicating poor performance under the evaluated criteria. DOGE, ETH, XLM, and LINK exhibit the highest sensitivity to methodological choices, with DOGE presenting the widest ranking range (positions 4--7).}
On the other hand, the correlation analysis in \autoref{fig:code_plot}(b) reveals important patterns in methodological agreement. Sum WSM and Vec WSM show \revone{Spearman} correlation of 1.000, indicating that these scaling methods produce identical rankings when paired with WSM aggregation. Similarly, Vec TOPSIS and Sum TOPSIS demonstrate very high correlation (0.983), suggesting robust agreement across different scaling approaches when using TOPSIS. However, the lowest correlations occur between WSM and TOPSIS methods paired with MinMax scaling (0.850), indicating that this represents the frontier of methodological sensitivity. Nevertheless, all correlations exceed 0.850, suggesting that while methodological decisions matter, the fundamental ranking structure remains relatively stable across this methodological space.
Finally, by directly analyzing the values of all generated rankings and their corresponding summary statistics in \autoref{tab:ranks}, the alternatives can be reasonably grouped into three categories according to their ranking variability:

\revone{
\begin{itemize}
\item \textbf{Perfect stability (MAD $= 0$):} Six alternatives exhibit zero median absolute deviation: ADA and XRP are completely fixed at positions 8 and 9 respectively; BNB and BTC alternate between the top two positions (1--2) depending on the scaling method; ETH maintains a median rank of 3, with occasional shifts to position 5 under specific configurations; and LTC occupies positions 6--7 with a median rank of 7.

\item \textbf{Moderate variability (MAD $= 0.5$):} XLM and LINK display moderate sensitivity to methodological choices. XLM fluctuates between positions 3 and 5 (median: 4), while LINK varies between positions 4 and 6 (median: 5.5).

\item \textbf{High variability (MAD $= 1.0$):} DOGE exhibits the highest sensitivity to methodological choices, presenting the widest ranking range (positions 4--7, median: 5.5).
\end{itemize}}

\section{Practical Implications}

Methodological variability becomes particularly critical when the evaluated alternatives compete near the decision boundary, as is often the case in funding selection processes.
Its relevance further increases in contexts where the costs of error are asymmetric, i.e., when an incorrect selection is more consequential than a justified rejection, and when transparent justification is required for diverse stakeholders.
Conversely, this factor is less influential when alternatives are clearly dominated or when decisions are preliminary in nature and allow for subsequent adjustments.

With respect to the use of combinatorial pipeline techniques, their application is especially advisable in high-impact scenarios such as strategic investments, public policy design, or budget allocation, where decision accuracy is essential.
In such contexts, the additional transparency provided by systematic methodological exploration is valuable for audits, scientific publications, and procurement processes.
These techniques are also particularly useful in exploratory phases where no clear methodological consensus exists, or in benchmarking studies aimed at comparing decision-making methods.

On the other hand, their application is not recommended for operational decisions requiring immediate responses (e.g., within less than one hour), nor in environments with limited hardware resources that preclude effective parallelization.
They should also be avoided when stakeholders have already agreed upon a well-justified preferred method, or when the methodological space becomes excessively large, exceeding, for instance, 100 combinations without a clear rationale to support such complexity.

\subsection{Usage Guidelines}

To implement combinatorial pipelines, we recommend starting with an initial exploration of the methodological space, applying \revone{Spearman rank correlation} analysis to identify redundancies.
In this stage, \revone{Spearman} correlation coefficients ($\rho>0.95$) are useful for detecting methods with nearly equivalent \revone{monotonic} behavior, while the coefficient of determination provides a deeper assessment by quantifying how accurately one method predicts another ($R^2>0.90$), or conversely, whether they capture distinct dimensions of the problem ($R^2<0.70$).
It is essential to quantify sensitivity through \revone{robust} summary statistics \revone{such as the median and the Median Absolute Deviation (MAD), which are more appropriate than the mean and standard deviation for ordinal ranking data}.

Positional stability can also be monitored using the \revone{normalized} Hamming distance; in this context, values of $d_H<0.20$ indicate that fewer than 20\% of alternatives change position, suggesting a stable ranking structure, whereas values of $d_H>0.40$ imply that more than 40\% of the list experiences positional shifts, representing substantive changes that warrant further examination.
Ultimately, the final method selection should be justified based on theory, existing literature, and interpretability, with ranking \revone{sensitivity} always reported alongside the final results.
It should be noted that these thresholds are illustrative, and analysts are expected to define their own critical limits according to the rigor and context of their study.

Finally, the exponential growth in the number of pipelines may render the process computationally demanding in large methodological spaces, motivating the use of intelligent sampling strategies as an alternative.
Moreover, given that the interpretation of statistical metrics requires substantial technical expertise, the use of visualization tools is essential for communicating results to non-expert users.
It should also be noted that the analysis relies on an implicit assumption of independence, whereby all methodological combinations are considered valid a priori; in practice, however, some combinations may be theoretically inappropriate or infeasible due to stakeholder constraints.
Lastly, since the proposed technique is validated within a single application domain, its generalization to other sectors requires further empirical validation, as different domains may exhibit distinct sensitivity patterns.

\section{Conclusions}
\label{section:conclusions}

In this work we present \texttt{SKCCombinatorialPipeline}, a framework that systematically addresses the problem of methodological \revtwo{sensitivity} in MCDM through automated exploration of the study space.
The implementation intelligently leverages the existing \textit{Scikit-Criteria} infrastructure, particularly \texttt{Ranks-Compara\allowbreak tor}, to provide robust comparative analysis of multiple methodological configurations.
The main contribution could lie in materializing theoretical principles of sensitivity analysis into practical tools that may facilitate exhaustive exploration of the methodological space without manual intervention, potentially establishing new standards of methodological transparency for contemporary MCDM practice.
A natural extension of this framework would be implementing adaptive sampling in large combinatorial spaces ($>100$ pipelines), centering the exploration in regions of high methodological variability.

\vspace{1pt}

\section*{Acknowledgments}

This work was supported by CONICET and CONAE.
A.R.S. was supported by CONICET fellowships.
This article has been revised using large language models (Claude,
ChatGPT, Gemini) in order to improve the clarity and correctness of the text.
All technical-scientific content remains property of the authors.
\normalsize

%
%
\bibliographystyle{splncs04}
\bibliography{bibliography}

\end{document}